\documentclass[draft,preprint,12pt,numbers,sort&compress]{elsarticle}
\usepackage{mathrsfs}
\usepackage{amssymb}
\usepackage{amsthm}
\usepackage{mathrsfs}
\usepackage[centertags]{amsmath}
\usepackage{amsfonts}

\usepackage{color}

\input amssym.def
\input amssym.tex

\date{}
\newtheorem{Definition}{Definition}[section]
\newtheorem{Theorem}{Theorem}[section]

\newtheorem{Lemma}{Lemma}[section]

\newcommand\R{\mbox{\bf R}}

\newcommand\SR{\mbox{\scriptsize\bf R}}

\newcommand{\definition}{{\lower .5ex
  \hbox{$\>\>\stackrel{\triangle}{=}\>\>$} }}

\begin{document}

\baselineskip=22pt
\thispagestyle{empty}
\baselineskip=22pt
\thispagestyle{empty}

\mbox{}
\bigskip

\begin{center}{\Large\bf Random data Cauchy problem for the wave}\\[1ex]
{\Large\bf   equation on   compact manifold}\\[4ex]
{Jinqiao  Duan$^a$,\quad Jianhua  Huang$^b$,Yongsheng Li$^c$,\quad Wei Yan\footnote{ Email: yanwei19821115@sina.cn}$^{d,a}$}\\[1ex]
{$^a$ Department of Applied Mathematics, Illinois Institute of Technology,}\\[1ex]
{Chicago, IL 60616, USA  }\\[1ex]

{$^b$College of Science, National University of Defense  Technology,}\\
{ Changsha, Hunan 410073,  China}\\[1ex]

{$^c$School of Mathematics, South  China  University of  Technology,}\\
{Guangzhou, Guangdong 510640,  China}\\[1ex]

{$^d$College of Mathematics and Information Science, Henan Normal University,}\\
{Xinxiang, Henan 453007,   China}\\[1ex]

\end{center}

\bigskip
\bigskip

\noindent{\bf Abstract.}
 Inspired by   the work of Burq and Tzvetkov
 (Invent. math. 173(2008), 449-475.), firstly,
 we construct the local strong solution to  the cubic nonlinear wave equation
with random data
for a large set of initial data  in $H^{s}(M)$ with $s\geq \frac{5}{14}$, where
 $M$ is a three dimensional compact  manifold with boundary, moreover, our result improves the result of Theorem 2 in (Invent. math. 173(2008), 449-475.);
  secondly,
we construct the local strong solution to
the quintic nonlinear wave equation
with random data for a large set of initial data  in $H^{s}(M)$ with $s\geq\frac{1}{6}$,
where $M$ is a two dimensional compact boundaryless  manifold;
finally,
we construct the local strong solution
to the  quintic nonlinear wave equation
with random data  for a large set of initial data  in $H^{s}(M)$ with $s\geq \frac{23}{90}$,
where $M$ is a two dimensional compact   manifold with  boundary.

\bigskip
\bigskip

{\large\bf 1. Introduction}
\bigskip

\setcounter{Theorem}{0} \setcounter{Lemma}{0}

\setcounter{section}{1}
In this paper, we investigate the  Cauchy problem for
\begin{eqnarray}
&& u_{tt}-\Delta u+u^{3}=0,\label{1.01}\\
&&(u,u_{t})\mid_{t=0}=(f_{1},f_{2}),\label{1.02}
\end{eqnarray}
with  real  initial  data  $f=(f_{1},f_{2})\in \mathscr{H}^{s}(M)=H^{s}(M)\times  H^{s-1}(M)$,
where $M$ is a three dimensional compact manifold with boundary.
We also investigate  the  Cauchy problem for
\begin{eqnarray}
&& u_{tt}-\Delta u+u^{5}=0,\label{1.03}\\
&&(u,u_{t})\mid_{t=0}=(f_{1},f_{2}),\label{1.04}
\end{eqnarray}
with  real  initial  data  $f=(f_{1},f_{2})\in \mathscr{H}^{s}(M)=H^{s}(M)\times  H^{s-1}(M)$,
where $M$ is a two dimensional compact manifold.

Rauch \cite{R} established the global regularity of (\ref{1.03})-(\ref{1.04}) in three dimension space with small initial energies.
Struwe \cite{S1998} obtained a unique global radially symmetric solution for any radially symmetric initial data $f_{1}\in C^{3}(\R^{3}),
f_{2}\in  C^{2}(\R^{3}).$ Some people  studied the global well-posedness, scattering and global space-time bounds \cite{BG,BS,GSV,G,GCPAM,K,KF,MP,N,N1999,Pe,Se,SS1994,T}.
Mockenhaupt et al. \cite{MSS} studied
local smoothing of  Fourier integral operators and Carleson-Sj\"olin estimates of the wave equation.
 Smith and Sogge \cite{SS1995} proved that   some Strichartz estimates for (\ref{1.03})-(\ref{1.04})
  hold on $n$ dimensional  Riemannian manifolds  with smooth, strictly geodesically concave
 boundaries and $n\geq2$. Christ et al. \cite{CCT}  proved that the solution map of (\ref{1.03})-(\ref{1.04})
 fails to be continuous at zero in the $H^{s}(\R^{3})\times H^{s-1}(\R^{3})$-topology for $0<s<1.$
 Christ et al. \cite{CCT} also proved that the solution map of (\ref{1.03})-(\ref{1.04})
 fails to be continuous at zero in the $H^{s}(\R^{2})\times H^{s-1}(\R^{2})$-topology for $0<s<\frac{1}{2}.$
 Burq  et al. \cite{BLP} studied   global existence for energy
critical waves in 3-D domains.

Burq and Tzvetkov \cite{BT2007} investigated the invariant measure for a three dimensional wave equation
\begin{eqnarray}
w_{tt}-\Delta w+|w|^{\alpha}w=0, (w,w_{t})|_{t=0}=(g_{1},g_{2}),\alpha<2\label{1.05}
\end{eqnarray}
with Dirichlet boundary condition and random initial data.
Burq and Tzvetkov \cite{BT2008-L} studied the local theory of the Cauchy problem for
(\ref{1.01})-(\ref{1.02}) with random data
in supercritical case on the three compact manifold. More precisely, they constructed
 the local strong solution
for a large set of initial data in $\mathscr{H}^{s}(M)$ with $s\geq\frac{1}{4}$, where $M$ is
 a three dimensional boundaryless compact manifold and
constructed the local strong solution
for a large set of initial data in $H^{s}(M)$ with $s\geq\frac{8}{21}$, where $M$ is a
 three dimensional  compact manifold with boundary; they also established the ill-posedness
 in $\mathscr{H}^{s}(M)$ with $s<\frac{1}{2}$ in the sense that
 the flow map on $\mathscr{H}^{s}(M)$ with $s<\frac{1}{2}$ is discontinuous at zero.
Burq and Tzvetkov \cite{BT2008-G}
obtained  the global existence result of the Cauchy problem for  a supercritical wave equation.
Bourgain and Bulut \cite{BBCR,BB,BB2014} studied Gibbs measure evolution in radial nonlinear wave on a three dimensional  ball.
Burq and Tzvetkov \cite{BT2014} established the probabilistic well-posedness for (\ref{1.01})
 in $H^{s}(M)$, $0<s<\frac{1}{2}$ with a suitable randomization on the three  dimensional torus.
L\"uhrmann and Mendelson \cite{LM} established an almost   sure    global   existence result of
 the defocusing nonlinear  wave equation of power-type on $\R^{3}$
with respect to a suitable randomization of the initial data. Recently,
Pocovnicu \cite{P} studied the   almost surely global well-posedness for the energy-critical defocusing
nonlinear wave equation on $\R^d, d=4$ and $5$ with random data.
 Very recently, Oh and  Pocovnicu \cite{OP} studied the  probabilistic global well-posedness of (\ref{1.01}) on ${\rm \R^{3}}$.

In this paper, inspired by \cite{BT2008-L,BSS,S}, firstly, for a large set of initial data  in $H^{s}(M)$ with $s\geq \frac{5}{14}$,
 we construct the local strong solution to the random data Cauchy problem for
(\ref{1.01}),
 where
 $M$ is a three dimensional compact  manifold with boundary;
  secondly,  for a large set of initial data  in $H^{s}(M)$ with $s\geq\frac{1}{6}$,
we construct the local strong solution
to the random data Cauchy problem for
 (\ref{1.03}),
where $M$ is a two dimensional compact boundaryless  manifold;
finally, for a large set of initial data  in $H^{s}(M)$ with $s\geq \frac{23}{90}$,
we construct the local strong solution
to the random data Cauchy problem for
 (\ref{1.03}),
  where $M$ is a two dimensional compact   manifold with  boundary.

We give some notations and some definitions before presenting the main results.
We assume that  $\Delta$ is the Laplace-Beltrami operator on  compact  manifold
and   $\Delta_{D}$ is the Laplace-Beltrami
 operator associated to the Dirichlet boundary condition
  and $\Delta_{N}$ is the Laplace-Beltrami
 operator associated to the Neumann  boundary condition.
   $(\Omega,\mathcal{F}, P)$  is a  probability space.   We define
 \begin{eqnarray*}
 \|v(t,x)\|_{L_{t}^{q}([0,T])L_{x}^{r}(M)}
 =\left(\int_{0}^{T}\left(\int_{M}|u|^{r}dx\right)^{\frac{q}{r}}dt\right)^{\frac{1}{q}}
 \end{eqnarray*}
and
 \begin{eqnarray*}
 \|u(\omega,t,x)\|_{L_{\omega}^{p}(\Omega)L_{t}^{q}([0,T])L_{x}^{r}(M)}
 =\left(\int_{\Omega}\left(\int_{0}^{T}\left(\int_{M}|u|^{r}dx\right)
 ^{\frac{q}{r}}dt\right)^{\frac{p}{q}}dP(\omega)\right)^{\frac{1}{p}}.
 \end{eqnarray*}

\begin{Definition}\label{1.1}
Assume that $(e_{n})\in C^{\infty}(M)(n=1,2,\cdot\cdot\cdot)$ is an orthonormal basis of  $L^{2}(M)$
and $(h_{n}(\omega),l_{n}(\omega))_{n=1}^{\infty}$ is a sequence of independent,
0 mean, real random variables on a probability space $(\Omega,\mathcal{F}, P)$ such that
\begin{eqnarray}
\exists C>0, \forall n\geq1, \int_{\Omega}\left(|h_{n}(\omega)|^{6}+|l_{n}(\omega)|^{6}\right)dP(\omega)<C.\label{1.06}
\end{eqnarray}
Let $f=(f_{1},f_{2})$, where
\begin{eqnarray*}
f_{1}=\sum_{n=1}^{\infty}\alpha_{n}e_{n}(x), f_{2}(x)
=\sum_{n=1}^{\infty}\beta_{n}e_{n}(x),\alpha_{n},\beta_{n}\in \R
\end{eqnarray*}
and the map
\begin{eqnarray}
\omega\longmapsto f^{\omega}=\left(f_{1}^{\omega}(x)
=\sum_{n=1}^{\infty}h_{n}(\omega)\alpha_{n}e_{n}(x),
 f_{2}(x)=\sum_{n=1}^{\infty}l_{n}(\omega)\beta_{n}e_{n}(x)\right)\label{1.07}
\end{eqnarray}
is equipped with the Borel sigma algebra from $(\Omega,\mathcal{F}).$
 From (\ref{1.07}), we know that the map $\omega\longmapsto f^{\omega}$
is measurable and $f^{\omega}\in L^{2}(\Omega; \mathscr{H}^{s}(M))$. Hence,
 this defines a $\mathscr{H}^{s}(M)$ valued random variable, which is the random function
related to $f$.
\end{Definition}

\begin{Definition}\label{1.2}
Assume that $M$  is  a smooth compact manifold and
$(e_{n})\in C^{\infty}(M)(n=1,2,\cdot\cdot\cdot)$ is an orthonormal basis of  $L^{2}(M)$
and $-\Delta e_{n}=\lambda_{n}^{2}e_{n}$. Let
\begin{eqnarray*}
H^{s}(M)=\left\{h\in H^{s}(M),h=\sum_{n=1}^{\infty}\gamma_{n}e_{n}(x),\|h\|_{H^{s}(M)}^{2}=
\sum_{n=1}^{\infty}(1+\lambda_{n}^{2})^{2s}|\gamma_{n}|^{2}<\infty\right\}.
\end{eqnarray*}
Define $\mathscr{H}^{s}(M)=H^{s}(M)\times H^{s-1}(M)$.
\end{Definition}
\begin{Definition}\label{1.3}
Assume that $M$  is  a smooth  manifold  with   boundary,  compact   closure and
$(e_{n})\in C^{\infty}(M)(n=1,2,\cdot\cdot\cdot)$ is an orthonormal basis of $L^{2}(M)$
and $-\Delta_{D} e_{n}=\lambda_{n}^{2}e_{n}$ with $e_{n}(x)\mid_{\partial M}=0.$ Let
\begin{eqnarray*}
H_{D}^{s}(M)=\left\{h=\sum_{n=1}^{\infty}\gamma_{n}e_{n}(x),\|h\|_{H_{D}^{s}(M)}^{2}=
\sum_{n=1}^{\infty}(1+\lambda_{n}^{2})^{2s}|\gamma_{n}|^{2}<\infty\right\}.
\end{eqnarray*}
Define $\mathscr{H}_{D}^{s}(M)=H_{D}^{s}(M)\times H_{D}^{s-1}(M)$.
\end{Definition}

\begin{Definition}\label{1.4}
Assume that $M$  is  a smooth  manifold  with   boundary,  compact   closure and
$(e_{n})\in C^{\infty}(M)(n=1,2,\cdot\cdot\cdot)$ is an orthonormal basis  of   $L^{2}(M)$
and $-\Delta_{N} e_{n}=\lambda_{n}^{2}e_{n}$ with $N_{x}. \nabla_{x}e_{n}(x)=0$,
 where $x\in \partial M$ and $N_{x}$ is a unit filed with respect to the metric. Let
\begin{eqnarray*}
H_{N}^{s}(M)=\left\{h=\sum_{n=1}^{\infty}\gamma_{n}e_{n}(x),\|h\|_{H_{N}^{s}(M)}^{2}
=\sum_{n=1}^{\infty}(1+\lambda_{n}^{2})^{2s}|\gamma_{n}|^{2}<\infty\right\}.
\end{eqnarray*}
Define $\mathscr{H}_{N}^{s}(M)=H_{N}^{s}(M)\times H_{N}^{s-1}(M)$.
\end{Definition}

The main result of this paper are as follows.
\begin{Theorem}\label{Thm1}
Let (\ref{1.06})   be  valid  and  $M$ be a three dimensional manifold with   boundary
 and $s\geq\frac{5}{14}$ and $f=(f_{1},f_{2})\in \mathscr{H}_{D}^{s}(M)$ and
 $f^{\omega}\in L^{2}(\Omega;\mathscr{H}_{D}^{s}(M)$ be defined by the randomization (\ref{1.07}).
 For a.s. $\omega \in \Omega$, there exist $T_{\omega}>0$ and a unique solution to (\ref{1.01}) with $u\mid _{R_{t}\times \partial M} =0$ and the
 initial data
 $f^{\omega}$ in a space continuously embedded in
 \begin{eqnarray*}
X_{\omega}=\left(cos(t\sqrt{-\Delta_{D}})f_{1}^{\omega}+\frac{sin (t\sqrt{-\Delta_{D}})}{\sqrt{-\Delta_{D}}}f_{2}^{\omega}\right)+
C([-T_{\omega}, T_{\omega}]; H^{\frac{2}{3}}_{D}(M)).
 \end{eqnarray*}
More precisely, for $0< T\leq1$, there exists $C>0,\delta>0$, an event $\Omega_{T}$ satisfying
\begin{eqnarray}
P(\Omega_{T})\geq 1-CT^{\frac{25}{14}}\label{1.08}
\end{eqnarray}
such that for every $\omega \in \Omega_{T}$ there exists a
unique solution of (\ref{1.01}) with data $f^{\omega}$ in a space continuously
embedded in $C([0,T]; H^{s}(M))$. Moreover, when  $h_{n},g_{n}$ are
standard real Gaussian or Bernoulli variables, we have
\begin{eqnarray}
P(\Omega_{T})\geq 1-C{\rm exp}\left(cT^{-\frac{1}{9}}\right).\label{1.09}
\end{eqnarray}
\end{Theorem}
\noindent {\bf Remark 1:} In Theorem 1.1,  if  $\Delta_{D}$,  $\mathscr{H}_{D}^{s}(M)$
  and Dirichlet boundary condition  $u\mid _{R_{t}\times \partial M} =0$
are replaced by
$\Delta_{N}$,  $\mathscr{H}_{N}^{s}(M)$ and  Neumann boundary condition
$N_{x}. \nabla_{x}u(x)\mid _{R_{t}\times \partial M} =0$, respectively, the conclusion is still  valid.
 $\mathscr{H}^{\frac{1}{2}}(M)$
is the critical space of (\ref{1.01})-(\ref{1.02}). In Theorem of \cite{BT2008-L}, the authors have constructed the local strong solution to
the cubic nonlinear wave equation with random data for a large set of initial data in $H^{s}(M)$ with $s\geq \frac{8}{21}.$
Thus, our result improves the result of \cite{BT2008-L}.

\begin{Theorem}\label{Thm2}
Let (\ref{1.06})   be  valid  and  $M$ be a two dimensional boundaryless manifold
 and $s\geq\frac{1}{6}$ and $f=(f_{1},f_{2})\in \mathscr{H}^{s}(M)$
  and $f^{\omega}\in L^{2}(\Omega;\mathscr{H}^{s}(M))$ be defined by
   the randomization (\ref{1.07}). For a.s. $\omega \in \Omega$, there exist
    $T_{\omega}>0$, $\sigma\geq \frac{1}{2}$ and a unique solution to (\ref{1.03}) with initial data
 $f^{\omega}$ in a space continuously embedded in
 \begin{eqnarray*}
X_{\omega}=\left(cos(t\sqrt{-\Delta})f_{1}^{\omega}
+\frac{sin (t\sqrt{-\Delta})}{\sqrt{-\Delta}}f_{2}^{\omega}\right)+
C([-T_{\omega}, T_{\omega}]; H^{\sigma}(M)).
 \end{eqnarray*}
More precisely, for $0< T\leq1$, there exists $C>0,\delta_{1}>0$, an event $\Omega_{T}$ satisfying
\begin{eqnarray}
P(\Omega_{T})\geq 1-CT^{1+\delta_{1}}\label{1.010}
\end{eqnarray}
 such that for every $\omega \in \Omega_{T}$ there exists a
 unique solution of (\ref{1.01}) with data $f^{\omega}$ in a space continuously
embedded in $C([0,T]; H^{s}(M))$. Moreover, when  $h_{n},g_{n}$
 are standard real Gaussian or Bernoulli variables, we have
\begin{eqnarray}
P(\Omega_{T})\geq 1-C{\rm exp}\left(cT^{-\delta_{1}}\right).\label{1.011}
\end{eqnarray}
\end{Theorem}

\begin{Theorem}\label{Thm3}
Let (\ref{1.06})   be  valid  and  $M$ be a three dimensional manifold with boundary
 and $s\geq\frac{23}{90}$ and $f=(f_{1},f_{2})\in \mathscr{H}_{D}^{s}(M)$ and
  $f^{\omega}\in L^{2}(\Omega;\mathscr{H}_{D}^{s}(M))$
  be defined by the randomization (\ref{1.07}). For a.s. $\omega \in \Omega$,
   there exist $T_{\omega}>0$ and a unique solution to (\ref{1.03}) with $u\mid _{R_{t}\times \partial M} =0$ and  the initial data
 $f^{\omega}$ in a space continuously embedded in
 \begin{eqnarray*}
X_{\omega}=\left(cos(t\sqrt{-\Delta_{D}})f_{1}^{\omega}+\frac{sin
(t\sqrt{-\Delta_{D}})}{\sqrt{-\Delta_{D}}}f_{2}^{\omega}\right)+
C([-T_{\omega}, T_{\omega}]; H^{\frac{7}{12}}_{D}(M)).
 \end{eqnarray*}
More precisely, for $0<T\leq1$, there exists $C>0,\delta_{2}>0$,
 an event $\Omega_{T}$ satisfying
\begin{eqnarray}
P(\Omega_{T})\geq 1-CT^{1+\delta_{2}}\label{1.012}
\end{eqnarray}
such that for every $\omega \in \Omega_{T}$ there exists a unique
 solution of (\ref{1.01})
with data $f^{\omega}$ in a space continuously
embedded in $C([0,T]; H^{s}(M))$. Moreover, when  $h_{n},g_{n}$ are
 standard real
 Gaussian or Bernoulli variables, we have
\begin{eqnarray}
P(\Omega_{T})\geq 1-C{\rm exp}\left(cT^{-\delta_{2}}\right).\label{1.013}
\end{eqnarray}
\end{Theorem}
\noindent {\bf Remark 2:} In Theorem 1.3,  if  $\Delta_{D}$,
  $\mathscr{H}_{D}^{s}(M)$
  and Dirichlet boundary condition $u\mid _{R_{t}\times \partial M} =0$
are replaced by
$\Delta_{N}$,  $\mathscr{H}_{N}^{s}(M)$ and
 Neumann boundary condition  $N_{x}. \nabla_{x}u(x)\mid _{R_{t}\times \partial M} =0$,
 respectively, the conclusion is still  valid.

The rest of the paper is arranged as follows. In Section 2,
 we give Strichartz
 estimates and $L^{p}(p=5,6)$ norm of eigenfunction associated
  to $-\Delta$ on compact manifold. In Section 3,
  we give some  properties of two random series. In Section 4,  we give averaging effects.
   In Section 5, we prove
the Theorem 1.1. In Section 6, we prove
the Theorem 1.2. In Section 7, we prove
the Theorem 1.3.

\bigskip
\bigskip

 \noindent{\large\bf 2.  Strichartz
 estimates and $L^{p}$ norm of eigenfunction associated  to $-\Delta$ on compact manifolds}
\setcounter{Definition}{0}
\setcounter{equation}{0}
\setcounter{Theorem}{0}
\setcounter{Lemma}{0}
\setcounter{section}{2}

In this section, we give some Strichartz
 estimates and $L^{p}$ norm of eigenfunction associated  to Laplace-Beltrami operator on
  two and three compact manifolds, which play a crucial role in establishing  Lemmas 4.1-4.6.

\begin{Definition}\label{2.001}
Let $0\leq s<1$ and $M$ be a three dimensional compact manifold with boundary.
A couple of real numbers $(p,q)$ is called $s$-admissible provided that
$p,q, s$ satisfy
\begin{eqnarray*}
\frac{1}{p}+\frac{3}{q}=\frac{3}{2}-s
\end{eqnarray*}
and $p\geq\frac{7}{2s}$ if $s\leq \frac{7}{10}$; $p=5$ if $s\geq \frac{7}{10}.$
For $T>0$, $0\leq s<1$, we define $X_{T}^{s}$ space and $Y_{T}^{s}$ space as follows.
\begin{eqnarray*}
&&X_{T}^{s}=C^{0}([0,T];H^{s}(M))\bigcap_{(p,q)s-admissible}L^{p}((0,T);L^{q}(M)),\\
&&Y_{T}^{s}=L^{1}([0,T];H^{-s}(M))+_{(p,q)s-admissible}L^{p^{\prime}}((0,T);L^{q^{\prime}}(M)),
\end{eqnarray*}
where $\frac{1}{p}+\frac{1}{p^{\prime}}=\frac{1}{q}+\frac{1}{q^{\prime}}=1.$
 Obviously, $Y_{T}^{s}$ is the dual space of $X_{T}^{s}.$
\end{Definition}

Definition 2.1 can be found in Definition 6.3 of \cite{BT2008-L}.

Inspired by (1.1)-(1.3) of \cite{BSS}, we give the definition 2.2.
\begin{Definition}\label{2.002}
Let $0\leq s<1$ and $M$ be a two dimensional compact   boundaryless manifold.
A couple of real numbers $(p,q)$ is called $s$-admissible provided that
$p,q, s$ satisfy
\begin{eqnarray*}
\frac{1}{p}+\frac{2}{q}=1-s
\end{eqnarray*}
and $\frac{3}{s}\leq p\leq \infty$.
For $T>0$, $0\leq s\leq1$, we define $X_{T}^{s}$ space and $Y_{T}^{s}$ space as follows.
\begin{eqnarray*}
&&X_{T}^{s}=C^{0}([0,T];H^{s}(M))\bigcap_{(p,q)s-admissible}L^{p}((0,T);L^{q}(M)),\\
&&Y_{T}^{s}=L^{1}([0,T];H^{-s}(M))+_{(p,q)s-admissible}L^{p^{\prime}}((0,T);L^{q^{\prime}}(M)),
\end{eqnarray*}
where $\frac{1}{p}+\frac{1}{p^{\prime}}=\frac{1}{q}+\frac{1}{q^{\prime}}=1.$ Obviously, $Y_{T}^{s}$ is the dual space of $X_{T}^{s}.$
\end{Definition}
Inspired by (1.1)-(1.2), (1.4) of \cite{BSS}, we give the definition 2.3.
\begin{Definition}\label{2.003}
Let $0\leq s<1$ and $M$ be a two dimensional compact manifold with boundary.
A couple of real numbers $(p,q)$ is called $s$-admissible provided that
$p,q, s$ satisfy
\begin{eqnarray*}
\frac{1}{p}+\frac{2}{q}=1-s
\end{eqnarray*}
and $p\geq\frac{5}{s}$ if $s\leq \frac{5}{8}$; $p=8$ if $s\geq \frac{5}{8}.$
For $T>0$, $0\leq s<1$, we define $X_{T}^{s}$ space and $Y_{T}^{s}$ space as follows.
\begin{eqnarray*}
&&X_{T}^{s}=C^{0}([0,T];H^{s}(M))\bigcap_{(p,q)s-admissible}L^{p}((0,T);L^{q}(M)),\\
&&Y_{T}^{s}=L^{1}([0,T];H^{-s}(M))+_{(p,q)s-admissible}L^{p^{\prime}}((0,T);L^{q^{\prime}}(M))
\end{eqnarray*}
where $\frac{1}{p}+\frac{1}{p^{\prime}}=\frac{1}{q}+\frac{1}{q^{\prime}}=1.$
 Obviously, $Y_{T}^{s}$ is the dual space of $X_{T}^{s}.$
\end{Definition}

\begin{Lemma}\label{Lemma2.1}
Let  $(p,q)$ be  an $s$-admissible couple of  Definition 2.1 and $M$ be a three dimensional compact manifold with boundary.
 For $0\leq s<1,$  there exists $C>0$ such that
\begin{eqnarray*}
&&\left\|cos (t\sqrt{-\Delta})(f_{1})\right\|_{X_{T}^{s}}
+\left\|\frac{sin(t\sqrt{-\Delta}) }{\sqrt{-\Delta}}(f_{2})
\right\|_{X_{T}^{s}}\leq C\|f\|_{\mathscr{H}^{s}(M)},\\
&&\left\|\int_{0}^{t}\frac{sin((t-\tau)\sqrt{-\Delta}) }
{\sqrt{-\Delta}}g(\tau)d\tau\right\|_{X_{T}^{s}}
\leq C\|g\|_{Y_{T}^{1-s}}
\end{eqnarray*}
for all $T\in (0,1]$ and $g\in \mathscr{H}^{s}(M).$
\end{Lemma}

Lemma 2.1 can be found in \cite{BT2008-L}.

\begin{Lemma}\label{Lemma2.2}
Let $(p,q)$ be an $s$-admissible couple of Definition 2.2 and $M$ be a
 two dimensional compact boundaryless  manifold.
   For $0\leq s<1,$ there exists $C>0$ such that
\begin{eqnarray*}
\left\|e^{\pm it\sqrt{-\Delta}}(g)\right\|_{L^{p}((0,T);L^{q}(M))}\leq C\|g\|_{H^{s}(M)}
\end{eqnarray*}
for all $T\in (0,1]$ and $g\in H^{s}(M).$
\end{Lemma}

For Lemma 2.2, we refer the readers to \cite{BSS}.

\begin{Lemma}\label{Lemma2.3}
Let  $(p,q)$ be  an $s$-admissible couple of  Definition 2.2 and $M$ be
 a two dimensional compact boundaryless manifold.
 For $0\leq s<1,$  there exists $C>0$ such that
\begin{eqnarray*}
&&\left\|cos (t\sqrt{-\Delta})(f_{1})\right\|_{X_{T}^{s}}
+\left\|\frac{sin(t\sqrt{-\Delta}) }{\sqrt{-\Delta}}(f_{2})
\right\|_{X_{T}^{s}}\leq C\|f\|_{\mathscr{H}^{s}(M)},\\
&&\left\|\int_{0}^{t}\frac{sin((t-\tau)\sqrt{-\Delta}) }
{\sqrt{-\Delta}}g(\tau)d\tau\right\|_{X_{T}^{s}}
\leq C\|g\|_{Y_{T}^{1-s}}
\end{eqnarray*}
for all $T\in (0,1]$ and $g\in \mathscr{H}^{s}(M).$
\end{Lemma}

Combining Lemma 2.2 with  the Corollary 4.3 of \cite{BT2007}, we have Lemma 2.3.

\begin{Lemma}\label{Lemma2.4}
Let $(p,q)$ be an $s$-admissible couple of Definition 2.3 and $M$ be a
 two dimensional compact   manifold with boundary.
   For $0\leq s<1,$ there exists $C>0$ such that
\begin{eqnarray*}
\left\|e^{\pm it\sqrt{-\Delta}}(f)\right\|_{L^{p}((0,T);L^{q}(M))}\leq C\|f\|_{H^{s}(M)}
\end{eqnarray*}
for all $T\in (0,1]$ and $f\in H^{s}(M).$
\end{Lemma}

For Lemma 2.4, we refer the readers to Theorem 1.1 of  \cite{BSS}.

\begin{Lemma}\label{Lemma2.5}
Let  $(p,q)$ be  an $s$-admissible couple of  Definition 2.3 and $M$ be
 a two dimensional compact  manifold with boundary.
 For $0\leq s<1,$  there exists $C>0$ such that
\begin{eqnarray*}
&&\left\|cos (t\sqrt{-\Delta})(f_{1})\right\|_{X_{T}^{s}}
+\left\|\frac{sin(t\sqrt{-\Delta}) }{\sqrt{-\Delta}}(f_{2})
\right\|_{X_{T}^{s}}\leq C\|f\|_{\mathscr{H}^{s}(M)},\\
&&\left\|\int_{0}^{t}\frac{sin((t-\tau)\sqrt{-\Delta}) }
{\sqrt{-\Delta}}g(\tau)d\tau\right\|_{X_{T}^{s}}
\leq C\|g\|_{Y_{T}^{1-s}}
\end{eqnarray*}
for all $T\in (0,1]$ and $g\in \mathscr{H}^{s}(M).$
\end{Lemma}

Combining Lemma 2.4 with  the Corollary of \cite{BT2007}, we have Lemma 2.5.

\begin{Lemma}\label{Lemma2.6}
Let $M$ be a three dimensional compact manifold with boundary and $(e_{n})_{n=1}^{\infty}$
 be an $L^{2}$-normalized basis consisting in eigenfunctions of the
 Laplace-Beltrami  operator with   Dirichlet (resp. Neumann)  boundary conditions, associated to eigenvalues
 $\lambda_{n}^{2}.$ Then,   there exists $C>0$ such that
\begin{eqnarray*}
\|e_{n}\|_{L^{5}(M)}\leq C(1+\lambda_{n}^{2})^{\frac{1}{5}}.
\end{eqnarray*}
\end{Lemma}

For the proof of Lemma 2.6, we refer the readers to Theorem 2 of \cite{SS}.

\begin{Lemma}\label{Lemma2.7}
Let $M$ be a two dimensional compact  boundaryless manifold and $(e_{n})_{n=1}^{\infty}$
 be an $L^{2}$-normalized basis consisting in eigenfunctions of the
 Laplace-Beltrami  operator, associated to eigenvalues
 $\lambda_{n}^{2}.$ Then,   there exists $C>0$ such that
\begin{eqnarray*}
\|e_{n}\|_{L^{6}(M)}\leq C(1+\lambda_{n}^{2})^{\frac{1}{12}}.
\end{eqnarray*}
\end{Lemma}

For the proof of Lemma 2.7, we refer the readers to Theorem 2.1 of \cite{S}.

\begin{Lemma}\label{Lemma2.8}
Let $M$ be a two dimensional compact  manifold with boundary and $(e_{n})_{n=1}^{\infty}$
 be an $L^{2}$-normalized basis consisting in eigenfunctions of the
 Laplace-Beltrami  operator with   Dirichlet (resp. Neumann)  boundary conditions, associated to eigenvalues
 $\lambda_{n}^{2}.$ Then,   there exists $C>0$ such that
\begin{eqnarray*}
\|e_{n}\|_{L^{6}(M)}\leq C(1+\lambda_{n}^{2})^{\frac{1}{9}}.
\end{eqnarray*}
\end{Lemma}

For the proof of Lemma 2.8, we refer the readers to Theorem 1.1 of \cite{SS}.

\bigskip
\bigskip

\noindent{\large\bf 3. Properties of two random series}

\setcounter{equation}{0}

 \setcounter{Theorem}{0}

\setcounter{Lemma}{0}

 \setcounter{section}{3}
In this section, we present $L^{p}$ properties of two random series which play a crucial role in proving
Lemmas 4.1-4.6.

 \begin{Lemma}\label{Lemma3.1}
Let $(l_{n}(\omega))_{n=1}^{\infty}$ be a sequence of independent,
 0-mean value,  complex random variables
 satisfying
\begin{eqnarray*}
\exists C>0, \forall n\geq 1,\left|\int_{\SR}|l_{n}(\omega)|^{2k}
dp(\omega)\right|\leq C.
\end{eqnarray*}
Then, we have that
\begin{eqnarray*}
&&\forall 2\leq p\leq 2k, \exists C>0, \forall (c_{n})_{n\in N^{*}}
\in l^{2}(N^{*},\mathbb{C}),\\
&&\left\|\sum_{n=1}^{\infty}c_{n}l_{n}\right\|_{L^{p}(\Omega)}\leq
 C\left(\sum_{n=1}^{\infty}|c_{n}|^{2}\right)^{\frac{1}{2}}.
\end{eqnarray*}
\end{Lemma}

Lemma 3.1 can be found in Lemma 4.2 of \cite{BT2008-L}.

\begin{Lemma}\label{Lemma3.2}
Let $(l_{n}(\omega))_{n=1}^{\infty}$ be a sequence of real, 0-mean, independent random variables with associated
sequence of distributions $(\mu_{n})_{n=1}^{\infty}$.
 Suppose that $\mu_{n}$ satisfy
\begin{eqnarray}
\exists C>0: \qquad \forall \gamma \in\ R, \forall n\geq 1,
\left|\int_{\SR}e^{\gamma x}d\mu_{n}(x)\right|\leq e^{C\gamma ^{2}}.\label{3.01}
\end{eqnarray}
Then, there exists $\alpha >0$ such that for every $\lambda >0$,
 every sequence $(c_{n})_{n=1}^{\infty}\in l^{2}$ of real
numbers,
\begin{eqnarray*}
P\left(\omega:\left|\sum_{n=1}^{\infty}c_{n}l_{n}(\omega)
\right|\right)\leq 2e^{-\frac{\alpha\lambda ^{2}}{\sum\limits_{n=1}^{\infty}c_{n}^{2}}}.
\end{eqnarray*}
Consequently, there exists $C>0$ such that
\begin{eqnarray*}
\left\|\sum_{n=1}^{\infty}c_{n}l_{n}(\omega)\right\|_{L^{p}(\Omega)}\leq
 \sqrt{p}\left(\sum_{n=1}^{\infty}c_{n}^{2}\right)^{\frac{1}{2}}
\end{eqnarray*}
for every $p\geq 2$ and every $(c_{n})_{n=1}^{\infty}\in l^{2}.$
\end{Lemma}

Lemma 3.2 can be found in Lemma 3.1 of \cite{BT2008-L}.

\bigskip
\bigskip

\noindent {\large\bf 4. Averaging effects}

\setcounter{equation}{0}

 \setcounter{Theorem}{0}

\setcounter{Lemma}{0}

\setcounter{section}{4}
In this section, motivated  by  Propositions 4.1, 4.4, 6.4 of \cite{BT2008-L},
we use Lemmas 2.6-2.8,  3.1, 3.2 to establish some mixed norm estimates about $u_{f}^{\omega}(x,t)$ defined below.

\begin{Lemma}\label{Lemma4.1}
Let $s\in \R$, $1<p\leq 5$ and $0<T\leq 1$ and  $f=(f_{1},f_{2})\in \mathscr{H}^{s}(M)$. Under the assumptions of Theorem 1, we have that
\begin{eqnarray}
\left\|(-\Delta+1)^{\frac{s}{2}-\frac{1}{5}}u_{f}^{\omega}\right\|_{L_{\omega}^{5}(\Omega)L_{t}^{p}([0,T])L_{x}^{5}(M)}\leq CT^{\frac{1}{p}}\|f\|_{\mathscr{H}^{s}(M)},\label{4.01}
\end{eqnarray}
where
$
u_{f}^{\omega}(x,t)=cos (t\sqrt{-\Delta})f_{1}^{\omega}+\frac{sin (t\sqrt{-\Delta})}{\sqrt{-\Delta}}f_{2}^{\omega}.
$
In particular,  for $s\in \R$, the following inequality is valid
\begin{eqnarray}
P\left(E_{\lambda, T,f}\right)\leq CT^{\frac{5}{p}}\lambda^{-5}\|f\|_{\mathscr{H}^{s}(M)}^5,\label{4.02}
\end{eqnarray}
where
$
E_{\lambda, T,f}=\left\{\omega\in \Omega:
\left\|(-\Delta+1)^{\frac{s}{2}-\frac{1}{8}}
u_{f}^{\omega}\right\|_{L_{t}^{p}([0,T])L_{x}^{5}(M)}\geq \lambda\right\}.
$
\end{Lemma}
\noindent{\bf Proof.} By using Lemma 3.1 and Minkowski inequality and Lemma 2.6, we have that
\begin{eqnarray}
&&\left\|cos(t\sqrt{-\Delta})f_{1}^{\omega}\right\|_{L_{\omega}^{5}
(\Omega)L_{t}^{p}([0,T])L_{x}^{5}(M)}\nonumber\\&&\leq C
\left\|\left(\sum_{n=1}^{\infty}|cos(t\lambda_{n})\alpha_{n}e_{n}(x)|^{2}
\right)^{\frac{1}{2}}\right\|_{L_{t}^{p}([0,T])L_{x}^{5}(M)}\nonumber\\
&&=\left\|\left(\sum_{n=1}^{\infty}|\alpha_{n}e_{n}(x)|^{2}\right)
^{\frac{1}{2}}\right\|_{L_{t}^{p}([0,T])L_{x}^{5}(M)}\nonumber\\
&&\leq CT^{\frac{1}{p}}\left[\left\|\sum_{n=1}^{\infty}|\alpha_{n}
e_{n}(x)|^{2}\right\|_{L_{x}^{\frac{5}{2}}(M)}^{\frac{1}{2}}\right]\nonumber\\
&&\leq CT^{\frac{1}{p}}\left[\sum_{n=1}^{\infty}|\alpha_{n}|^{2}
\|e_{n}(x)\|_{L^{5}}^{2}\right]^{\frac{1}{2}}\nonumber\\
&&\leq CT^{\frac{1}{p}}\left[\sum_{n=1}^{\infty}|\alpha_{n}|^{2}
(1+\lambda_{n}^{2})^{\frac{2}{5}}\right]^{\frac{1}{2}}
=CT^{\frac{1}{p}}\|f_{1}\|_{H^{\frac{2}{5}}(M)}.\label{4.03}
\end{eqnarray}
From (\ref{4.03}), for  $s\in\R$, we have that
\begin{eqnarray}
\left\|(-\Delta+1)^{\frac{s}{2}-\frac{1}{5}} cos(t\sqrt{-\Delta})f_{1}^{\omega}
\right\|_{L_{\omega}^{5}(\Omega)L_{t}^{p}([0,T])L_{x}^{5}(M)}\leq
 CT^{\frac{1}{p}}\|f_{1}\|_{H^{s}(M)}.\label{4.04}
\end{eqnarray}
From the property of  $\lambda_{n}^{2}(1\leq n\leq \infty,n\in N)$, we know that there exists $k\in N^{+}$ such that
\begin{eqnarray}
\lambda_{n}^{2}\leq 1,(1\leq n\leq k, n\in N);\lambda_{n}^{2}\geq 1, (n\geq k+1, n\in N).\label{4.05}
\end{eqnarray}
From (\ref{4.05}), for $0\leq t<1$ and $n\in N$,  we have that
\begin{eqnarray}
\left|\frac{sin t\lambda_{n}}{\lambda_{n}}\right|=|t|\left|\frac{sin t\lambda_{n}}{t\lambda_{n}}\right|\leq t\leq 1,(1\leq n\leq k);
\left|\frac{sin t\lambda_{n}}{\lambda_{n}}\right|\leq |\lambda_{n}|^{-1}, (n\geq k+1).\label{4.06}
\end{eqnarray}
 By using Lemma 3.1 and Minkowski inequality and Lemma 2.6 as well as (\ref{4.06}), we have that
\begin{eqnarray}
&&\left\|\frac{sin(t\sqrt{-\Delta})}{\sqrt{-\Delta}}f_{2}^{\omega}\right\|_{L_{\omega}^{5}
(\Omega)L_{t}^{p}([0,T])L_{x}^{5}(M)}\leq C
\left\|\left(\sum_{n=1}^{\infty}\left|\frac{sin(t\lambda_{n})}{\lambda_{n}}\beta_{n}e_{n}(x)\right|^{2}
\right)^{\frac{1}{2}}\right\|_{L_{t}^{p}([0,T])L_{x}^{5}(M)}\nonumber\\
&&\leq C\left\|\left\|\sum_{n=1}^{\infty}\left|\frac{sin(t\lambda_{n})}{\lambda_{n}}\beta_{n}
e_{n}(x)\right|^{2}\right\|_{L_{x}^{\frac{5}{2}}(M)}^{\frac{1}{2}}\right\|_{L_{t}^{p}([0,T])}\nonumber\\
&&\leq C\left\|\left[\sum_{n=1}^{\infty}\left|\frac{sin(t\lambda_{n})}{\lambda_{n}}\beta_{n}\right|^{2}
\|e_{n}(x)\|_{L^{5}}^{2}\right]^{\frac{1}{2}}\right\|_{L_{t}^{p}([0,T])}\nonumber\\
&&\leq C\left\|\left[\sum_{n=1}^{\infty}\left|\frac{sin(t\lambda_{n})}{\lambda_{n}}\beta_{n}\right|^{2}
(1+\lambda_{n}^{2})^{\frac{2}{5}}\right]^{\frac{1}{2}}\right\|_{L_{t}^{p}([0,T])}\nonumber\\
&&\leq\left\|\left[\sum_{n=1}^{k}\left|\frac{sin(t\lambda_{n})}{\lambda_{n}}\beta_{n}\right|^{2}
(1+\lambda_{n}^{2})^{\frac{2}{5}}\right]^{\frac{1}{2}}\right\|_{L_{t}^{p}([0,T])}+
\left\|\left[\sum_{n=k+1}^{\infty}\left|\frac{sin(t\lambda_{n})}{\lambda_{n}}\beta_{n}\right|^{2}
(1+\lambda_{n}^{2})^{\frac{2}{5}}\right]^{\frac{1}{2}}\right\|_{L_{t}^{p}([0,T])}\nonumber\\
&&\leq C\left\|\left[\sum_{n=1}^{k}\left|\beta_{n}\right|^{2}
(1+\lambda_{n}^{2})^{\frac{2}{5}}\right]^{\frac{1}{2}}\right\|_{L_{t}^{p}([0,T])}+
C\left\|\left[\sum_{n=k+1}^{\infty}\left|\beta_{n}\right|^{2}
(1+\lambda_{n}^{2})^{-\frac{3}{5}}\right]^{\frac{1}{2}}\right\|_{L_{t}^{p}([0,T])}\nonumber\\&&
\leq C\left\|\left[\sum_{n=1}^{k}\left|\beta_{n}\right|^{2}
(1+\lambda_{n}^{2})^{-\frac{3}{5}}\right]^{\frac{1}{2}}\right\|_{L_{t}^{p}([0,T])}+
C\left\|\left[\sum_{n=k+1}^{\infty}\left|\beta_{n}\right|^{2}
(1+\lambda_{n}^{2})^{-\frac{3}{5}}\right]^{\frac{1}{2}}\right\|_{L_{t}^{p}([0,T])}\nonumber\\
&&\leq CT^{\frac{1}{p}}\left[\sum_{n=1}^{\infty}\left|\beta_{n}\right|^{2}
(1+\lambda_{n}^{2})^{-\frac{3}{5}}\right]^{\frac{1}{2}}
=CT^{\frac{1}{p}}\|f_{1}\|_{H^{-\frac{3}{5}}(M)}.\label{4.07}
\end{eqnarray}
From (\ref{4.07}), we have that
\begin{eqnarray}
\left\|(-\Delta+1)^{\frac{s}{2}-\frac{1}{5}}\frac{sin(t\sqrt{-\Delta})}{\sqrt{-\Delta}}f_{2}^{\omega}\right\|_{L_{\omega}^{5}
(\Omega)L_{t}^{p}([0,T])L_{x}^{5}(M)}\leq CT^{\frac{1}{p}}\|f_{2}\|_{H^{s-1}(M)}.\label{4.08}
\end{eqnarray}
Combining (\ref{4.03})  with  (\ref{4.08}), we have that (\ref{4.01}) is valid.
By using the Bienaym$\acute{e}$-Tchebichev inequality and (\ref{4.01}), we have that
\begin{eqnarray}
&&P\left(E_{\lambda, T,f}\right)\leq \int_{\Omega} \frac{\left\|(-\Delta+1)^{\frac{s}{2}-\frac{1}{8}}
u_{f}^{\omega}\right\|_{L_{t}^{p}([0,T])L_{x}^{5}(M)}^{5}}{\lambda^{5}}dP(\omega)\nonumber\\&&\leq CT^{\frac{5}{p}}\lambda^{-5}\left\|(-\Delta+1)^{\frac{s}{2}-\frac{1}{5}}u_{f}^{\omega}\right\|_{L_{\omega}^{5}
(\Omega)L_{t}^{p}([0,T])L_{x}^{5}(M)}^{5}\nonumber\\&&\leq CT^{\frac{5}{p}}\lambda^{-5}\|f\|_{\mathscr{H}^{s}(M)}^{5}.\label{4.09}
\end{eqnarray}

This ends the proof of  Lemma 4.1.

\noindent {\bf Remark 3:} Our result improves the result of Proposition 6.4 of \cite{BT2008-L}.

\begin{Lemma}\label{Lemma4.2}
Let $s\in \R$, $1<q\leq 5$, $p\geq5$ and $0<T\leq 1$ and $f=(f_{1},f_{2})\in \mathscr{H}^{s}(M)$. Under the assumptions of Theorem 1,
if (\ref{3.01}) is valid, then we have that
\begin{eqnarray}
\left\|(-\Delta+1)^{\frac{s}{2}-\frac{1}{5}}u_{f}^{\omega}\right\|_{L_{\omega}^{p}(\Omega)L_{t}^{q}([0,T])L_{x}^{5}(M)}\leq CT^{\frac{1}{q}}\sqrt{p}\|f\|_{\mathscr{H}^{s}(M)},\label{4.010}
\end{eqnarray}
where
$
u_{f}^{\omega}(x,t)=cos (t\sqrt{-\Delta})f_{1}^{\omega}+\frac{sin (t\sqrt{-\Delta})}{\sqrt{-\Delta}}f_{2}^{\omega}.
$
In particular, for $s\in \R$,  the following inequality is valid
\begin{eqnarray*}
P\left(E_{\lambda,T, f}\right)\leq C{\rm exp}\left(-c\frac{\lambda^{2}}{\|f\|_{\mathscr{H}^{s}(M)}^{2}}\right),
\end{eqnarray*}
where
$
E_{\lambda, T,f}=\left\{\omega\in \Omega:
\left\|(-\Delta+1)^{\frac{s}{2}-\frac{1}{5}}
u_{f}^{\omega}\right\|_{L_{t}^{q}([0,T])L_{x}^{5}(M)}\geq \lambda\right\}.
$
\end{Lemma}

Combining Lemmas 3.2,  4.1 with the method of Proposition 4.4 of \cite{BT2008-L}, we derive that Lemma 4.2 is valid.

\begin{Lemma}\label{Lemma4.3}
Let $s\in \R$, $1<q\leq6$,  $0<T\leq 1$ and  $f=(f_{1},f_{2})\in \mathscr{H}^{s}(M)$.
Under the assumptions of Theorem 2, we have that
\begin{eqnarray}
\left\|(-\Delta+1)^{\frac{s}{2}-\frac{1}{12}}u_{f}^{\omega}\right\|_{L_{\omega}^{6}(\Omega)L_{t}^{q}([0,T])L_{x}^{6}(M)}\leq CT^{\frac{1}{q}}\|f\|_{\mathscr{H}^{s}(M)},\label{4.011}
\end{eqnarray}
where
$
u_{f}^{\omega}(x,t)=cos (t\sqrt{-\Delta})f_{1}^{\omega}+\frac{sin (t\sqrt{-\Delta})}{\sqrt{-\Delta}}f_{2}^{\omega}.
$
In particular, for $s\in \R$,  the following inequality is valid
\begin{eqnarray}
P\left(E_{\lambda, T,f}\right)\leq CT^{\frac{6}{q}}\lambda^{-6}\|f\|_{\mathscr{H}^{s}(M)}^{6},\label{4.012}
\end{eqnarray}
where
$
E_{\lambda, T,f}=\left\{\omega\in \Omega:
\left\|(-\Delta+1)^{\frac{s}{2}-\frac{1}{12}}
u_{f}^{\omega}\right\|_{L_{t}^{q}([0,T])L_{x}^{6}(M)}\geq \lambda\right\}.
$
\end{Lemma}

Lemma 4.3 can be proved similarly to Lemma 4.1 with the aid of Lemma 2.7.
\begin{Lemma}\label{Lemma4.4}
Let $s\in \R$, $1<q\leq6,$  $p\geq5,$ and $0<T\leq 1$ and $f=(f_{1},f_{2})\in \mathscr{H}^{s}(M)$.
 Under the assumptions of Theorem 2,
  if (\ref{3.01}) is valid, then we have that
\begin{eqnarray}
\left\|(-\Delta+1)^{\frac{s}{2}-\frac{1}{12}}u_{f}^{\omega}
\right\|_{L_{\omega}^{p}(\Omega)L_{t}^{q}([0,T])L_{x}^{6}(M)}
\leq CT^{\frac{1}{q}}\|f\|_{\mathscr{H}^{s}(M)},\label{4.013}
\end{eqnarray}
where
$
u_{f}^{\omega}(x,t)=cos (t\sqrt{-\Delta})
f_{1}^{\omega}+\frac{sin (t\sqrt{-\Delta})}{\sqrt{-\Delta}}f_{2}^{\omega}.
$
In particular, for $s\in \R$,  the following inequality is valid
\begin{eqnarray*}
P\left(E_{\lambda,T, f}\right)\leq C{\rm exp}
\left(-c\frac{\lambda^{2}}{\|f\|_{\mathscr{H}^{s}(M)}^{2}}\right),
\end{eqnarray*}
where
$
E_{\lambda, T,f}=\left\{\omega\in \Omega:
\left\|(-\Delta+1)^{\frac{s}{2}-\frac{1}{12}}
u_{f}^{\omega}\right\|_{L_{t}^{q}([0,T])L_{x}^{6}(M)}\geq \lambda\right\}.
$
\end{Lemma}

Lemma 4.4 can be proved similarly to Lemma 4.2 with the aid of Lemma 4.3.

\begin{Lemma}\label{Lemma4.5}
Let $s\in \R$, $1<q\leq6$, $0<T\leq 1$ and  $f=(f_{1},f_{2})\in \mathscr{H}^{s}(M)$.
 Under the assumptions of Theorem 3, we have that
\begin{eqnarray}
\left\|(-\Delta+1)^{\frac{s}{2}-\frac{1}{12}}u_{f}^{\omega}
\right\|_{L_{\omega}^{6}(\Omega)L_{t}^{q}([0,T])
L_{x}^{6}(M)}\leq CT^{\frac{1}{q}}\|f\|_{\mathscr{H}^{s}(M)},\label{4.014}
\end{eqnarray}
where
$
u_{f}^{\omega}(x,t)=cos (t\sqrt{-\Delta})f_{1}^{\omega}
+\frac{sin (t\sqrt{-\Delta})}{\sqrt{-\Delta}}f_{2}^{\omega}.
$
In particular, for $s\in \R$,  the following inequality is valid
\begin{eqnarray}
P\left(E_{\lambda, T,f}\right)\leq CT^{\frac{6}{q}}
\lambda^{-6}\|f\|_{\mathscr{H}^{s}(M)}^{6},\label{4.015}
\end{eqnarray}
where
$
E_{\lambda, T,f}=\left\{\omega\in \Omega:
\left\|(-\Delta+1)^{\frac{s}{2}-\frac{1}{9}}
u_{f}^{\omega}\right\|_{L_{t}^{q}([0,T])L_{x}^{6}(M)}\geq \lambda\right\}.
$
\end{Lemma}

Lemma 4.5 can be proved similarly to Lemma 4.1 with the aid of Lemma 2.8.
\begin{Lemma}\label{Lemma4.6}
Let $s\in \R$, $1<q\leq6$ and $p\geq6,$ $0<T\leq 1$, $f=(f_{1},f_{2})\in \mathscr{H}^{s}(M)$.
 Under the assumptions of Theorem 3,  if (\ref{3.01}) is valid, then we have that
\begin{eqnarray}
\left\|(-\Delta+1)^{\frac{s}{2}-\frac{1}{9}}u_{f}^{\omega}
\right\|_{L_{\omega}^{p}(\Omega)L_{t}^{q}([0,T])L_{x}^{6}(M)}
\leq CT^{\frac{1}{q}}\|f\|_{\mathscr{H}^{s}(M)},\label{4.016}
\end{eqnarray}
where
$
u_{f}^{\omega}(x,t)=cos (t\sqrt{-\Delta})f_{1}^{\omega}
+\frac{sin (t\sqrt{-\Delta})}{\sqrt{-\Delta}}f_{2}^{\omega}.
$
In particular, for $s\in \R$,  the following inequality is valid
\begin{eqnarray*}
P\left(E_{\lambda,T, f}\right)\leq C{\rm exp}
\left(-c\frac{\lambda^{2}}{\|f\|_{\mathscr{H}^{s}(M)}^{2}}\right),
\end{eqnarray*}
where
$
E_{\lambda, T,f}=\left\{\omega\in \Omega:
\left\|(-\Delta+1)^{\frac{s}{2}-\frac{1}{9}}
u_{f}^{\omega}\right\|_{L_{t}^{q}([0,T])L_{x}^{6}(M)}\geq \lambda\right\}.
$
\end{Lemma}

Lemma 4.6 can be proved similarly to Lemma 4.2.

\bigskip
\bigskip

\noindent {\large\bf 5. Proof of Theorem 1.1}

\setcounter{equation}{0}

 \setcounter{Theorem}{0}

\setcounter{Lemma}{0}

\setcounter{section}{5}
In this section, following  the  method of  \cite{BT2008-L},  we use  Lemmas 4.1, 4.2 and  contraction map theorem to  prove Theorem 1.1.
We give Lemma 5.1 before proving Theorem 1.1.

\begin{Lemma}\label{Lemma5.1}
Let $s\geq \frac{5}{14}$ and $0<T\leq1$. Then, we have that
  \begin{eqnarray*}
 K_{f}^{\omega}(v)=-\int_{0}^{t}\frac{sin((t-\tau)\sqrt{-\Delta})}
 {\sqrt{-\Delta}}\left((u_{f}^{\omega}+v)^{3}\right)(\cdot,\tau)d\tau
  \end{eqnarray*}
  and $(v,v_{t})|_{t=0}=(0,0).$
  Then, there exists $C>0$  such that for every
   $f\in \mathscr{H}^{s}(M)$ and $\omega \in E_{\lambda,f}^{c}$,
    the map $K_{f}^{\omega}$ satisfies
  \begin{eqnarray}
  &&\|K_{f}^{\omega}(u)\|_{X_{T}^{\frac{2}{3}}}\leq
   C\left(\lambda^{3}+T^{\frac{1}{3}}\|v\|_{X_{T}^{\frac{1}{3}}}^{3}\right),\label{5.01}\\
  &&\|K_{f}^{\omega}(v)-K_{f}^{\omega}(w)\|_{X_{T}^{\frac{2}{3}}}
  \leq CT^{\frac{1}{9}}\left[\lambda^{2}+\|v\|_{X_{T}^{\frac{2}{3}}}^{2}
  +\|w\|_{X_{T}^{\frac{2}{3}}}^{2}\right]\label{5.02}.
  \end{eqnarray}
\end{Lemma}
\noindent{\bf Proof.}
Since $(\frac{21}{4},\frac{14}{3})$ is $\frac{2}{3}$-admissible, we have that
\begin{eqnarray}
\|g\|_{L^{\infty}([0,T]; H^{\frac{2}{3}}(M))}+\|g\|_{L_{t}^{\frac{21}{4}}
([0,T];L^{\frac{14}{3}}(M))}\leq C\|g\|_{X_{T}^{\frac{2}{3}}}.\label{5.03}
\end{eqnarray}
By using (\ref{5.03}) and Lemma 2.1, we have that
\begin{eqnarray}
&&\left\|K_{f}^{\omega}(u)\right\|_{X_{T}^{\frac{2}{3}}}\leq C\left\|\int_{0}^{t}\frac{sin((t-\tau)\sqrt{-\Delta})}{\sqrt{-\Delta}}\left((u_{f}^{\omega}+v)^{3}\right)
(\cdot,\tau)d\tau\right\|_{X_{T}^{\frac{2}{3}}}\nonumber\\&&\leq C\|\left((u_{f}^{\omega}+v)^{3}\right)\|_{Y_{T}^{\frac{1}{3}}}\leq C\left\|(u_{f}^{\omega}+v)^{3}\right\|_{L_{t}^{\frac{21}{19}}([0,T];L_{x}^{\frac{14}{9}}(M))}\nonumber\\
&&\leq C\left(\left\|u_{f}^{\omega}\right\|_{L_{t}^{\frac{63}{19}}([0,T];L_{x}^{\frac{14}{3}}(M))}^{3}
+\left\|v\right\|_{L_{t}^{\frac{63}{19}}([0,T];L_{x}^{\frac{14}{3}}(M))}^{3}\right)\label{5.04}.
\end{eqnarray}
By using  the H\"older inequality,  from (\ref{5.03}), we have that
\begin{eqnarray}
\left\|v\right\|_{L_{t}^{\frac{63}{19}}([0,T];L_{x}^{\frac{14}{3}}(M))}\leq C T^{\frac{1}{9}}\left\|v\right\|_{L_{t}^{\frac{21}{4}}([0,T];L_{x}^{\frac{14}{3}}(M))}
\leq CT^{\frac{1}{9}}\|v\|_{X_{T}^{\frac{2}{3}}}\label{5.05}.
\end{eqnarray}
When $s\geq1$, since $5(s-\frac{2}{5})\geq 3$, we have that
\begin{eqnarray}
W^{s-\frac{2}{5},5}(M)\hookrightarrow L^{\frac{14}{3}}(M)\label{5.06}.
\end{eqnarray}
By using (\ref{5.06}), for  $\omega \in E_{\lambda,T,f}^{c},$  we have that
\begin{eqnarray}
\left\|u_{f}^{\omega}\right\|_{L_{t}^{\frac{63}{19}}([0,T];L_{x}^{\frac{14}{3}}(M))}\leq C\left\|u_{f}^{\omega}\right\|_{L_{t}^{\frac{63}{19}}([0,T];W^{s-\frac{2}{5},5}(M))}\leq
 \lambda\label{5.07}.
\end{eqnarray}
When $s<1$, by using the Sobolev embedding Theorem, we have
\begin{eqnarray}
W^{s-\frac{2}{5},5}(M)\hookrightarrow L^{q_{1}}(M),
\frac{1}{q_{1}}=\frac{1}{5}-\frac{s-\frac{2}{5}}{3}=\frac{1-s}{3}.\label{5.08}
\end{eqnarray}
By using (\ref{5.08}), for $\omega \in E_{\lambda,T,f}^{c}$, we have that
\begin{eqnarray}
\left\|u_{f}^{\omega}\right\|_{L_{t}^{\frac{63}{19}}([0,T];L_{x}^{q_{1}}(M))}\leq C\left\|u_{f}^{\omega}\right\|_{L_{t}^{\frac{63}{19}}([0,T];W^{s-\frac{2}{5},5}(M))}\leq C\lambda\label{5.09}.
\end{eqnarray}
Since $\frac{14}{3}\leq \frac{3}{1-s}$ which is equivalent to $s\geq \frac{5}{14}$, from (\ref{5.09}), we have that
\begin{eqnarray}
\left\|u_{f}^{\omega}\right\|_{L_{t}^{\frac{63}{19}}([0,T];L_{x}^{\frac{14}{3}}(M))}\leq \left\|u_{f}^{\omega}\right\|_{L_{t}^{\frac{63}{19}}([0,T];L_{x}^{q_{1}}(M))}\leq C\left\|u_{f}^{\omega}\right\|_{L_{t}^{\frac{63}{19}}([0,T];W^{s-\frac{2}{5},5}(M))}\leq C\lambda\label{5.010}.
\end{eqnarray}
Inserting (\ref{5.05}), (\ref{5.07}), (\ref{5.09})-(\ref{5.010}) into (\ref{5.04}) yields that
\begin{eqnarray}
&&\left\|K_{f}^{\omega}(u)\right\|_{X_{T}^{\frac{2}{3}}}\leq C\left(\lambda^{3}
+T^{\frac{1}{3}}\left\|v\right\|_{X_{T}^{\frac{2}{3}}}^{3}\right)\label{5.011}.
\end{eqnarray}
By using a  proof similar to (\ref{5.011}),  we have that
\begin{eqnarray}
  \|K_{f}^{\omega}(v)-K_{f}^{\omega}(w)\|_{X_{T}^{\frac{2}{3}}}\leq
   CT^{\frac{1}{9}}\|v-w\|_{X_{T}^{\frac{2}{3}}}
  \left[\lambda^{2}+T^{\frac{2}{9}}\|v\|_{X_{T}^{\frac{2}{3}}}^{2}
  +T^{\frac{2}{9}}\|w\|_{X_{T}^{\frac{2}{3}}}^{2}\right]\label{5.012}.
 \end{eqnarray}

This completes the proof of Lemma 5.1.

Now we prove Theorem 1.1.
Fix $0<T\leq1$. Let
\begin{eqnarray}
B(0, 2C\lambda^{3})=\left\{u|u\in X_{T}^{\frac{2}{3}}:\|u\|_{ X_{T}^{\frac{2}{3}}}
\leq 2C\lambda^{3}\right\}\label{5.013}
\end{eqnarray}
and
\begin{eqnarray}
4CT^{\frac{1}{9}}\lambda^{3}\leq \lambda.\label{5.014}
\end{eqnarray}
From (\ref{5.011})-(\ref{5.014}), we have that
\begin{eqnarray}
&&\left\|K_{f}^{\omega}(u)\right\|_{X_{T}^{\frac{2}{3}}}\leq C\left(\lambda^{3}
+T^{\frac{1}{3}}\left\|v\right\|_{X_{T}^{\frac{2}{3}}}^{3}\right)\leq C\left(\lambda^{3}
+T^{\frac{1}{3}}\left(2C\lambda^{3}\right)^{3}\right)\leq 2C\lambda^{3},\label{5.015}\\
&&\|K_{f}^{\omega}(v)-K_{f}^{\omega}(w)\|_{X_{T}^{\frac{2}{3}}}\leq CT^{\frac{1}{9}}
\|v-w\|_{X_{T}^{\frac{2}{3}}}
  \left[\lambda^{2}+T^{\frac{2}{9}}\|v\|_{X_{T}^{\frac{2}{3}}}^{2}+T^{\frac{2}{9}}
  \|w\|_{X_{T}^{\frac{2}{3}}}^{2}\right]\nonumber\\
  &&\leq \frac{1}{2}\|v-w\|_{X_{T}^{\frac{2}{3}}}.\label{5.016}
 \end{eqnarray}
Thus, $K_{f}^{\omega}$ is a contraction map on the ball $B(0, 2C\lambda^{3})$.
We define
\begin{eqnarray}
\Omega_{T}=E^{c}_{\lambda, T,f}, \sum=\bigcup\limits_{n\in N^{*}}\Omega_{\frac{1}{n}}\label{5.017}.
\end{eqnarray}
Combining (\ref{4.02}) with (\ref{5.017}),  we have that
\begin{eqnarray}
P(\Omega_{T})\geq 1-CT^{\frac{25}{14}}, P\left(\sum\right)=1.\label{5.018}
\end{eqnarray}
When  $h_{n},g_{n}$ are standard real Gaussian or Bernoulli variables, (\ref{3.01}) is valid.
In this case, by using a proof similar to case (\ref{1.06}), we have that $K_{f}^{\omega}$
 is a contraction map on the ball $B(0, 2C\lambda^{3})$. We define
 \begin{eqnarray}
\Omega_{T}=E^{c}_{\lambda, T,f}, \sum=\bigcup\limits_{n\in N^{*}}\Omega_{\frac{1}{n}}.\label{5.019}
\end{eqnarray}
Combining (\ref{5.019}) with Lemma 4.2, we have that
\begin{eqnarray*}
P(\Omega_{T})\geq1-C{\rm exp}\left(-cT^{-\frac{1}{9}}\right), P\left(\sum\right)=1.
\end{eqnarray*}

This completes the proof of Theorem 1.1.

\bigskip
\bigskip

\noindent {\large\bf 6. Proof of Theorem 1.2}

\setcounter{equation}{0}

 \setcounter{Theorem}{0}

\setcounter{Lemma}{0}

\setcounter{section}{6}

In this section, following  the  method of  \cite{BT2008-L}, we use  Lemmas 4.3, 4.4 and
  contraction map theorem to  prove Theorem 1.2.
We present Lemmas 6.1, 6.2 before proving Theorem 1.2.

\begin{Lemma}\label{Lemma6.1}
Let $s=\frac{1}{6}$ and $0<T\leq1$. Then, we have that
  \begin{eqnarray*}
 K_{f}^{\omega}(v)=-\int_{0}^{t}
 \frac{sin((t-\tau)\sqrt{-\Delta})}{\sqrt{-\Delta}}\left((u_{f}^{\omega}+v)^{5}\right)(\cdot,\tau)d\tau
  \end{eqnarray*}
  and $(v,v_{t})|_{t=0}=(0,0).$
  Then, there exists $C>0$  such that for every $f\in \mathscr{H}^{s}(M)$
   and $\omega \in E_{\lambda,f}^{c}$, the map $K_{f}^{\omega}$ satisfies
  \begin{eqnarray}
  &&\|K_{f}^{\omega}(u)\|_{X_{T}^{\frac{1}{2}}}\leq C\left(\lambda^{5}
  +\|v\|_{X_{T}^{\frac{1}{2}}}^{5}\right),\label{6.01}\\
  &&\|K_{f}^{\omega}(v)-K_{f}^{\omega}(w)\|_{X_{T}^{\frac{1}{2}}}\leq C
  \|v-w\|_{X_{T}^{\frac{1}{2}}}\left[\lambda^{4}+\|v\|_{X_{T}^{\frac{1}{2}}}^{4}
  +\|w\|_{X_{T}^{\frac{1}{2}}}^{4}\right]\label{6.02}.
  \end{eqnarray}
\end{Lemma}
\noindent{\bf Proof.}For $\omega \in E_{\lambda,T,f}^{c}$, we have that
\begin{eqnarray}
\left\|u_{f}^{\omega}\right\|_{L^{6}((0,T)\times M)}\leq \lambda,\label{6.03}
\end{eqnarray}
from (\ref{6.03}) and  Lemma 2.3, we have that
\begin{eqnarray}
\left\|K_{f}^{\omega}(v)\right\|_{X_{T}^{\frac{1}{2}}}\leq C
\left\|\left(u_{f}^{\omega}+v\right)^{5}\right\|_{L^{\frac{6}{5}}([0,T]\times M)}\leq C
\left[\left\|u_{f}^{\omega}\right\|_{L^{6}([0,T]\times M)}^{5}
+\left\|v\right\|_{L^{6}([0,T]\times M)}^{5}\right]\label{6.04}.
\end{eqnarray}
Since $(6,6)$ is $\frac{1}{2}$-admissible, combining (\ref{6.03})  with  (\ref{6.04}),  we have that
\begin{eqnarray}
\left\|K_{f}^{\omega}(v)\right\|_{X_{T}^{\frac{1}{2}}}\leq C
\left\|\left(u_{f}^{\omega}+v\right)^{5}\right\|_{L^{\frac{6}{5}}([0,T]\times M)}\leq C
\left[\lambda^{5}+\left\|v\right\|_{X_{T}^{\frac{1}{2}}}^{5}\right]\label{6.05}.
\end{eqnarray}
By using a proof similar to (\ref{6.05}), we have that
\begin{eqnarray}
  \|K_{f}^{\omega}(v)-K_{f}^{\omega}(w)\|_{X_{T}^{\frac{1}{2}}}\leq C
  \|v-w\|_{X_{T}^{\frac{1}{2}}}\left[\lambda^{4}+\|v\|_{X_{T}^{\frac{1}{2}}}^{4}
  +\|w\|_{X_{T}^{\frac{1}{2}}}^{4}\right]\label{6.06}.
  \end{eqnarray}

This completes the proof of Lemma 6.1.

\begin{Lemma}\label{Lemma6.2}
Let $s\geq\frac{1}{6}+\frac{4\epsilon}{15}$,  $0<\epsilon\ll\frac{1}{2}$ and $0<T\leq1$. Then, we have that
  \begin{eqnarray*}
 K_{f}^{\omega}(v)=-\int_{0}^{t}\frac{sin((t-\tau)\sqrt{-\Delta})}{\sqrt{-\Delta}}
 \left((u_{f}^{\omega}+v)^{5}\right)(\cdot,\tau)d\tau
  \end{eqnarray*}
  and $(v,v_{t})|_{t=0}=(0,0).$
  Then, there exists $C>0$  such that for every $f\in \mathscr{H}^{s}(M)$
   and $\omega \in E_{\lambda,T,f}^{c}$,
  the map $K_{f}^{\omega}$ satisfies
  \begin{eqnarray}
  &&\|K_{f}^{\omega}(u)\|_{X_{T}^{\frac{1}{2}+\epsilon}}\leq C
  \left(\lambda^{5}+T^{\frac{22\epsilon}{3}}\|v\|_{X_{T}^{\frac{1}{2}+\epsilon}}^{5}\right),
  \label{6.07}\\
  &&\|K_{f}^{\omega}(v)-K_{f}^{\omega}(w)\|_{X_{T}^{\frac{1}{2}+\epsilon}}
  \leq CT^{\frac{22\epsilon}{15}}\left[\lambda^{4}+\|v_{1}\|_{X_{T}^{\frac{1}{2}+\epsilon}}^{4}
  +\|w\|_{X_{T}^{\frac{1}{2}+\epsilon}}^{4}\right]\label{6.08}.
  \end{eqnarray}
\end{Lemma}
\noindent{\bf Proof.}
Since $(\frac{30}{5-22\epsilon},\frac{30}{5-4\epsilon})$ is $(\frac{1}{2}+\epsilon)$-admissible, we have that
\begin{eqnarray}
\|g\|_{L^{\infty}([0,T]; H^{\frac{1}{2}+\epsilon}(M))}+\|g\|_{L_{t}^{\frac{30}{5-22\epsilon}}
([0,T];L^{\frac{30}{5-4\epsilon}}(M))}\leq C\|g\|_{X_{T}^{\frac{1}{2}+\epsilon}}.\label{6.09}
\end{eqnarray}
By using (\ref{6.09}), we have that
\begin{eqnarray}
&&\left\|K_{f}^{\omega}(u)\right\|_{X_{T}^{\frac{1}{2}+\epsilon}}\leq C\left\|\int_{0}^{t}\frac{sin((t-\tau)\sqrt{-\Delta})}{\sqrt{-\Delta}}\left((u_{f}^{\omega}+v)^{5}\right)
(\cdot,\tau)d\tau\right\|_{X_{T}^{\frac{1}{2}+\epsilon}}\nonumber\\&&\leq C
\|\left((u_{f}^{\omega}+v)^{3}\right)\|_{Y_{T}^{\frac{1}{2}-\epsilon}}\leq C\left\|(u_{f}^{\omega}+v)^{5}\right\|_{L_{t}^{\frac{6}{5+2\epsilon}}([0,T];L_{x}^{\frac{6}{5-4\epsilon}}(M))}\nonumber\\
&&\leq C\left(\left\|u_{f}^{\omega}\right\|_{L_{t}^{\frac{30}{5+2\epsilon}}([0,T];L_{x}^{\frac{30}{5-4\epsilon}}(M))}^{5}
+\left\|v\right\|_{L_{t}^{\frac{30}{5+2\epsilon}}([0,T];L_{x}^{\frac{30}{5-4\epsilon}}(M))}^{5}\right)\label{6.010}.
\end{eqnarray}
By using  the H\"older inequality,  from (\ref{6.09}), we have that
\begin{eqnarray}
\left\|v\right\|_{L_{t}^{\frac{30}{5+22\epsilon}}([0,T];L_{x}^{\frac{30}{5-4\epsilon}}(M))}\leq C T^{\frac{22\epsilon}{15}}\left\|v\right\|_{L_{t}^{\frac{30}{5-22\epsilon}}([0,T];L_{x}^{\frac{30}{5-4\epsilon}}(M))}\leq CT^{\frac{22}{15}\epsilon}\|v\|_{X_{T}^{\frac{1}{2}+\epsilon}}\label{6.011}.
\end{eqnarray}
When $s\geq\frac{1}{2}$, since $6(s-\frac{1}{6})\geq 2$, we have that
\begin{eqnarray}
W^{s-\frac{1}{6},6}(M)\hookrightarrow L^{\frac{30}{5-4\epsilon}}(M)\label{6.012}.
\end{eqnarray}
By using (\ref{6.012}), for  $\omega \in E_{\lambda,T,f}^{c},$  we have that
\begin{eqnarray}
\left\|u_{f}^{\omega}\right\|_{L_{t}^{\frac{30}{5+22\epsilon}}([0,T];L_{x}^{\frac{30}{5-4\epsilon}}(M))}\leq C\left\|u_{f}^{\omega}\right\|_{L_{t}^{\frac{30}{5+22\epsilon}}([0,T];W^{s-\frac{2}{5},5}(M))}
\leq \lambda\label{6.013}.
\end{eqnarray}
When $s<\frac{1}{2}$, by using the Sobolev embedding Theorem, we have
\begin{eqnarray}
W^{s-\frac{1}{6},6}(M)\hookrightarrow L^{q_{1}}(M),\frac{1}{q_{1}}=\frac{1}{6}-\frac{s-\frac{1}{6}}{2}
=\frac{1-2s}{4}.\label{6.014}
\end{eqnarray}
By using (\ref{6.014}), for $\omega \in E_{\lambda,T,f}^{c}$, we have that
\begin{eqnarray}
\left\|u_{f}^{\omega}\right\|_{L_{t}^{\frac{30}{5+22\epsilon}}([0,T];L_{x}^{q_{1}}(M))}\leq C\left\|u_{f}^{\omega}\right\|_{L_{t}^{\frac{30}{5+22\epsilon}}([0,T];W^{s-\frac{1}{6},6}(M))}\leq C\lambda\label{6.015}.
\end{eqnarray}
Since $\frac{30}{5-4\epsilon}\leq \frac{2}{1-2s}$ which is equivalent to $s\geq \frac{5+8\epsilon}{30}$, from (\ref{6.015}), we have that
\begin{eqnarray}
&&\left\|u_{f}^{\omega}\right\|_{L_{t}^{\frac{30}{5+22\epsilon}}([0,T];L_{x}^{\frac{30}{5-4\epsilon}}(M))}\leq \left\|u_{f}^{\omega}\right\|_{L_{t}^{\frac{30}{5+22\epsilon}}([0,T];L_{x}^{q_{1}}(M))}\nonumber\\&&\leq C\left\|u_{f}^{\omega}\right\|_{L_{t}^{\frac{30}{5+22\epsilon}}([0,T];W^{s-\frac{2}{5},5}(M))}\leq C\lambda\label{6.016}.
\end{eqnarray}
Inserting (\ref{6.011}), (\ref{6.013}), (\ref{6.015})-(\ref{6.016}) into (\ref{6.010}) yields that
\begin{eqnarray}
&&\|K_{f}^{\omega}(u)\|_{X_{T}^{\frac{1}{2}+\epsilon}}\leq C
  \left(\lambda^{5}+T^{\frac{22\epsilon}{3}}\|v\|_{X_{T}^{\frac{1}{2}+\epsilon}}^{5}\right)\label{6.017}.
\end{eqnarray}
By using a  proof similar to (\ref{6.017}),  we have that
\begin{eqnarray}
  \|K_{f}^{\omega}(v)-K_{f}^{\omega}(w)\|_{X_{T}^{\frac{2}{3}}}\leq CT^{\frac{22\epsilon}{15}}
  \|v-w\|_{X_{T}^{\frac{1}{2}+\epsilon}}\left[\lambda^{4}
  +\|v\|_{X_{T}^{\frac{2}{3}}}^{4}+\|w\|_{X_{T}^{\frac{1}{2}+\epsilon}}^{4}\right]\label{6.018}.
 \end{eqnarray}

This completes the proof of Lemma 6.2.

Now we prove Theorem 1.2.
Fix $0<T\leq1$. Firstly, we consider $s=\frac{1}{6}.$ Let
\begin{eqnarray}
B(0, 2C\lambda^{5})=\left\{u|u\in X_{T}^{\frac{1}{2}}:\|u\|_{ X_{T}^{\frac{1}{2}}}\leq 2C\lambda^{5}\right\}\label{6.019}
\end{eqnarray}
and
\begin{eqnarray}
4C\lambda^{5}\leq \lambda.\label{6.020}
\end{eqnarray}
From (\ref{6.01})-(\ref{6.02}), we have that
\begin{eqnarray}
&&\left\|K_{f}^{\omega}(u)\right\|_{X_{T}^{\frac{1}{2}}}\leq C\left(\lambda^{5}
+\left\|v\right\|_{X_{T}^{\frac{1}{2}}}^{5}\right)\leq C\left(\lambda^{5}
+\left(2C\lambda^{5}\right)^{5}\right)\leq 2C\lambda^{5},\label{6.021}\\
&&\|K_{f}^{\omega}(v)-K_{f}^{\omega}(w)\|_{X_{T}^{\frac{1}{2}}}\leq C
\|v-w\|_{X_{T}^{\frac{1}{2}}}
  \left[\lambda^{4}+\|v\|_{X_{T}^{\frac{1}{2}}}^{4}+
  \|w\|_{X_{T}^{\frac{1}{2}}}^{4}\right]\nonumber\\
  &&\leq \frac{1}{2}\|v-w\|_{X_{T}^{\frac{1}{2}}}.\label{6.022}
 \end{eqnarray}
Thus, $K_{f}^{\omega}$ is a contraction map on the ball $B(0, 2C\lambda^{5})$.
We define
\begin{eqnarray}
\Omega_{T}=E^{c}_{\lambda, T,f}, \sum=\bigcup\limits_{n\in N^{*}}\Omega_{\frac{1}{n}}\label{6.023}.
\end{eqnarray}
Combining (\ref{4.02}) with (\ref{6.023}),  we have that
\begin{eqnarray}
P(\Omega_{T})\geq 1-CT, P\left(\sum\right)=1.\label{6.024}
\end{eqnarray}
We consider case $s>\frac{1}{6}.$
Fix $0<T\leq1$. Let
\begin{eqnarray}
B(0, 2C\lambda^{3})=\left\{u|u\in X_{T}^{\frac{1}{2}+\epsilon}:\|u\|_{ X_{T}^{\frac{1}{2}+\epsilon}}\leq
 2C\lambda^{5}\right\}\label{6.025}
\end{eqnarray}
and
\begin{eqnarray}
2CT^{\frac{22\epsilon}{15}}\lambda^{5}\leq \lambda.\label{6.026}
\end{eqnarray}
From (\ref{6.017})-(\ref{6.018}), we have that
\begin{eqnarray}
&&\left\|K_{f}^{\omega}(u)\right\|_{X_{T}^{\frac{1}{2}+\epsilon}}\leq C\left(\lambda^{5}
+T^{\frac{22\epsilon}{3}}\left\|v\right\|_{X_{T}^{\frac{1}{2}+\epsilon}}^{5}\right)\leq C\left(\lambda^{5}
+T^{\frac{22\epsilon}{3}}\left(2C\lambda^{5}\right)^{5}\right)\leq 2C\lambda^{5},\label{6.027}\\
&&\|K_{f}^{\omega}(v)-K_{f}^{\omega}(w)\|_{X_{T}^{\frac{1}{2}+\epsilon}}\leq CT^{\frac{88\epsilon}{15}}
\|v-w\|_{X_{T}^{\frac{1}{2}+\epsilon}}
  \left[\lambda^{4}+T^{\frac{88\epsilon}{15}}\|v\|_{X_{T}^{\frac{1}{2}+\epsilon}}^{4}+T^{\frac{88\epsilon}{15}}
  \|w\|_{X_{T}^{\frac{1}{2}+\epsilon}}^{4}\right]\nonumber\\
  &&\leq \frac{1}{2}\|v-w\|_{X_{T}^{\frac{1}{2}+\epsilon}}.\label{6.028}
 \end{eqnarray}
Thus, $K_{f}^{\omega}$ is a contraction map on the ball $B(0, 2C\lambda^{5})$.
We define
\begin{eqnarray}
\Omega_{T}=E^{c}_{\lambda, T,f}, \sum=\bigcup\limits_{n\in N^{*}}\Omega_{\frac{1}{n}}\label{6.029}.
\end{eqnarray}
Combining (\ref{4.02}) with (\ref{6.029}),  we have that
\begin{eqnarray}
P(\Omega_{T})\geq 1-CT^{1+\frac{33\epsilon}{5}}, P\left(\sum\right)=1.\label{6.030}
\end{eqnarray}
When  $h_{n},g_{n}$ are standard real Gaussian or Bernoulli variables, (\ref{3.01}) is valid.
In this case, by using a proof similar to case (\ref{1.06}), we have that $K_{f}^{\omega}$
 is a contraction map on the ball $B(0, 2C\lambda^{5})$. We define
 \begin{eqnarray}
\Omega_{T}=E^{c}_{\lambda, T,f}, \sum=\bigcup\limits_{n\in N^{*}}\Omega_{\frac{1}{n}}.\label{6.031}
\end{eqnarray}
Combining (\ref{6.031})   with   Lemma 4.4, we have that
\begin{eqnarray*}
P(\Omega_{T})\geq1-C{\rm exp}\left(-cT^{-\frac{11\epsilon}{15}}\right), P\left(\sum\right)=1.
\end{eqnarray*}

This completes the proof of  Theorem 1.2.

\bigskip
\bigskip
\noindent {\large\bf 7. Proof of Theorem 1.3}

\setcounter{equation}{0}

 \setcounter{Theorem}{0}

\setcounter{Lemma}{0}

\setcounter{section}{7}
In this section, following  the  method of  \cite{BT2008-L},  we use  Lemmas 4.5, 4.6 and  contraction map theorem to  prove Theorem 1.3.
We give Lemma 7.1 before proving Theorem 1.3.

\begin{Lemma}\label{Lemma7.1}
Let $s\geq \frac{23}{90}$ and $0<T\leq1$. Then, we have that
  \begin{eqnarray*}
 K_{f}^{\omega}(v)=-\int_{0}^{t}\frac{sin((t-\tau)\sqrt{-\Delta})}
 {\sqrt{-\Delta}}\left((u_{f}^{\omega}+v)^{5}\right)(\cdot,\tau)d\tau
  \end{eqnarray*}
  and $(v,v_{t})|_{t=0}=(0,0).$
  Then, there exists $C>0$  such that for every $f\in \mathscr{H}^{s}(M)$
   and $\omega \in E_{\lambda,T,f}^{c}$, the map $K_{f}^{\omega}$ satisfies
  \begin{eqnarray}
  &&\hspace{-0.8cm}\|K_{f}^{\omega}(u)\|_{X_{T}^{\frac{7}{12}}}\leq C\left(\lambda^{5}
  +T^{\frac{1}{3}}\|v\|_{X_{T}^{\frac{1}{3}}}^{5}\right),\label{7.01}\\
  &&\hspace{-0.8cm}\|K_{f}^{\omega}(v)-K_{f}^{\omega}(w)\|_{X_{T}^{\frac{7}{12}}}\leq
   CT^{\frac{1}{15}}\|v-w\|_{X_{T}^{\frac{7}{12}}}\left[\lambda^{4}
  +T^{\frac{4}{15}}\|v\|_{X_{T}^{\frac{7}{12}}}^{4}+T^{\frac{4}{15}}
  \|w\|_{X_{T}^{\frac{7}{12}}}^{4}\right]\label{7.02}.
  \end{eqnarray}
\end{Lemma}
\noindent{\bf Proof.}
Since $(\frac{60}{7},\frac{20}{3})$ is $\frac{7}{12}$-admissible, we have that
\begin{eqnarray}
\|g\|_{L^{\infty}([0,T]; H^{\frac{7}{12}}(M))}+\|g\|_{L_{t}^{\frac{60}{7}}
([0,T];L^{\frac{20}{3}}(M))}\leq C\|g\|_{X_{T}^{\frac{7}{12}}}.\label{7.03}
\end{eqnarray}
By using (\ref{7.03}) and Lemma 2.5, we have that
\begin{eqnarray}
&&\left\|K_{f}^{\omega}(u)\right\|_{X_{T}^{\frac{7}{12}}}\leq C
\left\|\int_{0}^{t}\frac{sin((t-\tau)\sqrt{-\Delta})}{\sqrt{-\Delta}}
\left((u_{f}^{\omega}+v)^{5}\right)
(\cdot,\tau)d\tau\right\|_{X_{T}^{\frac{7}{12}}}\nonumber\\&&\leq C
\|\left((u_{f}^{\omega}+v)^{5}\right)
\|_{Y_{T}^{\frac{5}{12}}}\leq C\left\|(u_{f}^{\omega}+v)^{5}\right\|_{L_{t}^{\frac{12}{11}}
([0,T];L_{x}^{\frac{4}{3}}(M))}
\nonumber\\
&&\leq C\left(\left\|u_{f}^{\omega}\right\|
_{L_{t}^{\frac{60}{11}}([0,T];L_{x}^{\frac{20}{3}}(M))}^{5}
+\left\|v\right\|_{L_{t}^{\frac{60}{11}}([0,T];L_{x}^{\frac{20}{3}}(M))}^{5}
\right)\label{7.04}.
\end{eqnarray}
By using  the H\"older inequality,  from (\ref{7.03}), we have that
\begin{eqnarray}
\left\|v\right\|_{L_{t}^{\frac{60}{11}}([0,T];L_{x}^{\frac{20}{3}}(M))}\leq C T^{\frac{1}{15}}\left\|v\right\|_{L_{t}^{\frac{60}{7}}([0,T];L_{x}^{\frac{20}{3}}(M))}
\leq CT^{\frac{1}{15}}\|v\|_{X_{T}^{\frac{7}{12}}}\label{7.05}.
\end{eqnarray}
When $s\geq\frac{5}{9}$, since $6(s-\frac{2}{9})\geq 2$, we have that
\begin{eqnarray}
W^{s-\frac{2}{9},6}(M)\hookrightarrow L^{\frac{14}{3}}(M)\label{7.06}.
\end{eqnarray}
By using (\ref{7.06}), for  $\omega \in E_{\lambda,T,f}^{c},$  we have that
\begin{eqnarray}
\left\|u_{f}^{\omega}\right\|_{L_{t}^{\frac{60}{11}}([0,T];L_{x}^{\frac{20}{3}}(M))}\leq C\left\|u_{f}^{\omega}\right\|_{L_{t}^{\frac{60}{11}}([0,T];W^{s-\frac{2}{9},6}(M))}\leq
\lambda\label{7.07}.
\end{eqnarray}
When $s<\frac{5}{9}$, by using the Sobolev embedding Theorem, we have that
\begin{eqnarray}
W^{s-\frac{2}{9},6}(M)\hookrightarrow L^{q_{1}}(M),\frac{1}{q_{1}}
=\frac{1}{6}-\frac{s-\frac{2}{9}}{2}=\frac{1-s}{3}.\label{7.08}
\end{eqnarray}
By using (\ref{7.08}), for $\omega \in E_{\lambda,T,f}^{c}$, we have that
\begin{eqnarray}
\left\|u_{f}^{\omega}\right\|_{L_{t}^{\frac{60}{11}}([0,T];L_{x}^{q_{1}}(M))}\leq C\left\|u_{f}^{\omega}\right\|_{L_{t}^{\frac{60}{11}}([0,T];W^{s-\frac{2}{9},6}(M))}
\leq C\lambda\label{7.09}.
\end{eqnarray}
Since $\frac{20}{3}\leq \frac{18}{5-9s}$ which is equivalent to $s\geq \frac{23}{90}$,
from (\ref{7.09}), we have that
\begin{eqnarray}
\left\|u_{f}^{\omega}\right\|_{L_{t}^{\frac{60}{11}}([0,T];L_{x}^{\frac{14}{3}}(M))}\leq \left\|u_{f}^{\omega}\right\|_{L_{t}^{\frac{60}{11}}([0,T];L_{x}^{q_{1}}(M))}\leq C\left\|u_{f}^{\omega}\right\|_{L_{t}^{\frac{60}{11}}([0,T];W^{s-\frac{2}{9},6}(M))}
\leq C\lambda\label{7.010}.
\end{eqnarray}
Inserting (\ref{7.05}), (\ref{7.07}), (\ref{7.09})-(\ref{7.010}) into (\ref{7.04}) yields that
\begin{eqnarray}
&&\left\|K_{f}^{\omega}(u)\right\|_{X_{T}^{\frac{7}{12}}}\leq C\left(\lambda^{5}
+T^{\frac{1}{5}}\left\|v\right\|_{X_{T}^{\frac{7}{12}}}^{5}\right)\label{7.011}.
\end{eqnarray}
By using a  proof similar to (\ref{7.011}),  we have that
\begin{eqnarray}
  \|K_{f}^{\omega}(v)-K_{f}^{\omega}(w)\|_{X_{T}^{\frac{7}{12}}}\leq CT^{\frac{1}{15}}
  \|v-w\|_{X_{T}^{\frac{2}{3}}}\left[\lambda^{4}
  +\|v\|_{X_{T}^{\frac{7}{12}}}^{4}+\|w\|_{X_{T}^{\frac{7}{12}}}^{4}\right]\label{7.012}.
 \end{eqnarray}

This completes the proof of Lemma 7.1.

Now we prove Theorem 1.3.

By using Lemmas 7.1, 4.5, 4.6 and a proof similar to Theorem 1.1, we can obtain Theorem 1.3.

\leftline{\large \bf Acknowledgments}

\bigskip

\noindent

  The fourth author is  supported by
  the Young core Teachers Program of Henan Normal University.

\end{document}